\documentclass[a4paper,11pt,twoside,reqno]{amsart}

\usepackage[utf8]{inputenc}
\usepackage[plainpages=false,pdfpagelabels=true]{hyperref}
\usepackage{amssymb,amsthm}
\usepackage[margin=1in]{geometry}

\newtheorem{Satz}{Theorem}[section]
\newtheorem{Prop}[Satz]{Proposition}
\newtheorem{Lem}[Satz]{Lemma}

\newtheorem{Cor}[Satz]{Corollary}
\theoremstyle{definition}
\newtheorem{Dfn}[Satz]{Definition}
\newtheorem{Bem}[Satz]{Remark}
\newtheorem{Bsp}[Satz]{Example}
\newcommand{\tr}{\operatorname{Tr}}

\newcommand{\C}{{\mathbb{C}}}
\newcommand{\vol}{{\operatorname{vol}}}
\newcommand{\dv}{\text{ }dV}
\allowdisplaybreaks[1]

\renewcommand{\epsilon}{\varepsilon}

\newcommand{\R}{\ensuremath{\mathbb{R}}}

\newcommand{\cL}{\mathcal{L}}

\numberwithin{equation}{section}

\title{On interpolating sesqui-harmonic maps between Riemannian manifolds}
\author{Volker Branding}
\date{\today}
\address{University of Vienna, Faculty of Mathematics\\
Oskar-Morgenstern-Platz 1, 1090 Vienna, Austria\\}
\email{volker.branding@univie.ac.at}
\begin{document}

\begin{abstract}
Motivated from the action functional for bosonic strings with
extrinsic curvature term we introduce an action functional for maps 
between Riemannian manifolds that interpolates between the actions
for harmonic and biharmonic maps.
Critical points of this functional will be called interpolating sesqui-harmonic maps.
In this article we initiate a rigorous mathematical treatment of this functional
and study various basic aspects of its critical points.
\end{abstract} 

\maketitle

\section{Introduction and Results}
\emph{Harmonic maps} play an important role in geometry, analysis and physics.
On the one hand they are one of the most studied variational problems in geometric analysis,
on the other hand they naturally appear in various branches of theoretical physics,
for example as critical points of the nonlinear sigma model or in the theory of elasticity.
Mathematically, they are defined as critical points of the Dirichlet energy
\begin{align}
\label{harmonic-energy}
E_1(\phi)=\int_M|d\phi|^2\dv,
\end{align}
where \(\phi\colon M\to N\) is a map between the two Riemannian manifolds \((M,h)\) and \((N,g)\).
The critical points of \eqref{harmonic-energy} are characterized by the vanishing of the so-called
\emph{tension field}, which is given by
\begin{align*}
0=\tau(\phi):=\tr_h\nabla d\phi.
\end{align*}
This is a semilinear, elliptic second order partial differential equation, for which many results 
on existence and qualitative behavior of its solutions have been obtained.
For a recent survey on harmonic maps see \cite{MR2389639}.
Due to their nonlinear nature harmonic maps do not always need to exist.
For example, if \(M=T^2\) and \(N=S^2\), there does not exist a harmonic map with \(\deg\phi=\pm 1\) 
regardless of the chosen metrics \cite{MR0420708}.

In these cases one may consider the following generalization of the harmonic map equation,
the so-called \emph{biharmonic maps}.
These arise as critical points of the bienergy \cite{MR886529}, which is defined as
\begin{align}
\label{biharmonic-energy}
E_2(\phi)=\int_M|\tau(\phi)|^2\dv.
\end{align}
In contrast to the harmonic map equation, the biharmonic map equation is an elliptic equation of fourth order 
and is characterized by the vanishing of the \emph{bitension field}
\begin{align*}
0=\tau_2(\phi):=\Delta\tau(\phi)-R^N(d\phi(e_\alpha),\tau(\phi))d\phi(e_\alpha),
\end{align*}
where \(\Delta\) is the connection Laplacian on \(\phi^\ast TN\), \(e_\alpha\) an orthonormal basis of \(TM\) and \(R^N\) denotes the curvature tensor
of the target manifold \(N\). We make use of the Einstein summation convention, meaning
that we sum over repeated indices. 

In the literature that studies analytical aspects of biharmonic maps one refers to \eqref{biharmonic-energy}
as the energy functional for \emph{intrinsic biharmonic maps}.

For a survey on biharmonic maps between Riemannian manifolds we refer to
\cite{MR3098705} and \cite{MR2301373}.

In this article we want to focus on the study of an action functional
that interpolates between the actions for harmonic and biharmonic maps

\begin{align}
\label{energy-functional}
E_{\delta_1,\delta_2}(\phi)=\delta_1\int_M|d\phi|^2\dv+\delta_2\int_M|\tau(\phi)|^2\dv
\end{align}
with \(\delta_1, \delta_2\in\R\).

This functional appears at several places in the physics literature.
In string theory it is known as \emph{bosonic string with extrinsic curvature term},
see \cite{MR848552,MR834521}.

On the mathematical side there have been several articles dealing
with some particular aspect of \eqref{energy-functional}.
Up to the best knowledge of the author the first place where the functional
\eqref{energy-functional} was mentioned is \cite[pp.134-135]{MR0440602}
with \(\delta_2=1\) and \(\delta_1>0\).
In that reference it is already shown that if the domain has dimension $2$ or $3$
and the target \(N\) negative sectional curvature then the critical points of 
\eqref{energy-functional} reduce to harmonic maps.
Later it was shown in \cite[p.191]{MR636281}
that no critical points exist if one does not impose the curvature condition
on \(N\) and also assumes that \(\deg\phi=1\).
Some analytic questions related to critical points of \eqref{energy-functional} have been discussed in \cite{MR2124627}
assuming \(\delta_1=2,\delta_2=1\). 
For the sake of completeness we want to mention that the functional
\eqref{energy-functional} with \(\delta_1>0\) and \(\delta_2=\frac{1}{2}\) is also 
presented in the survey article ``A report on harmonic maps'', 
see \cite[p.28, Example (6.30)]{MR495450}.

In \cite{MR2465916} the authors initiate an extensive study of
\eqref{energy-functional} assuming \(\delta_2=1\) and \(\delta_1\in\R\)
under the condition that \(\phi\) is an immersion. They consider variations of 
\eqref{energy-functional} that are normal to the image \(\phi(M)\subset N\).
In this setup they call critical points of \eqref{energy-functional}
\emph{biminimal immersions}. They also point out possible applications
of their model to the theory of elasticity.

Up to now there exist several results on biminimal immersions, see
for example \cite{MR2838514} for biminimal hypersurfaces into spheres,
\cite{MR3231742} for biminimal submanifolds in manifolds of non-positive curvature
and \cite{MR2964661} for biminimal submanifolds of Euclidean space.
Instead of investigating maps that are immersions, we here want to put the focus
on arbitrary maps between Riemannian manifolds.

The critical points of \eqref{energy-functional} will be referred to as \emph{interpolating sesqui-harmonic maps} and are given by
\begin{align}
\label{euler-lagrange}
\delta_2\Delta\tau(\phi)=\delta_2 R^N(d\phi(e_\alpha),\tau(\phi))d\phi(e_\alpha)+\delta_1\tau(\phi),
\end{align}
where \(\tau(\phi)\) denotes the tension field of the map \(\phi\)
and by \(\Delta\) we are representing
the connection Laplacian on the vector bundle \(\phi^\ast TN\).
\par\medskip
As in the case of biharmonic maps it is obvious that harmonic maps solve \eqref{euler-lagrange}.
For this reason we are mostly interested in solutions of \eqref{euler-lagrange} that are not harmonic maps.
However, we can expect that as in the case of biharmonic maps there may be 
many situations in which solutions of \eqref{euler-lagrange} will be harmonic maps.
In particular, we can expect that this is the case if \(N\) has negative sectional curvature
and \(\delta_1\delta_2>0\). This question will be dealt with in section 4.
On the other hand, if \(\delta_1\) and \(\delta_2\) have opposite sign we might
expect a different behavior of solutions of \eqref{euler-lagrange} 
since in this case the two terms in the energy functional \eqref{energy-functional}
are competing with each other and the energy functional can become unbounded 
from above and below.
\par\medskip
This article is organized as follows:
In section 2 we study basic features of interpolating sesqui-harmonic maps.
Afterwards, in section 3, we derive several explicit solutions of the interpolating sesqui-harmonic 
map equation and in the last section we provide several results that characterize
the qualitative behavior of interpolating sesqui-harmonic maps.

Throughout this paper we will make use of the following conventions.
Whenever choosing local coordinates we will use Greek letters to 
denote indices on the domain \(M\) and Latin letters for
indices on the target \(N\). We will choose the following convention
for the curvature tensor \(R(X,Y)Z:=[\nabla_X,\nabla_Y]Z-\nabla_{[X,Y]}Z\)
such that the sectional curvature is given by \(K(X,Y)=R(X,Y,Y,X)\).
For the Laplacian acting on functions \(f\in C^\infty(M)\) we choose the convention
\(\Delta f=\operatorname{div}\operatorname{grad}f\),
for sections in the vector bundle \(\phi^{\ast}TN\) we make
the choice \(\Delta^{\phi^\ast TN}=\tr(\nabla^{\phi^\ast TN}\nabla^{\phi^\ast TN})\).
Note that the connection Laplacian on \(\phi^\ast TN\) is defined by
\(\Delta:=\nabla_{e_\alpha}\nabla_{e_\alpha}-\nabla_{\nabla_{e_\alpha}e_\alpha}\).

\section{Interpolating sesqui-harmonic maps}
In this section we analyze the basic features of the action functional \eqref{energy-functional} 
and start by calculating its critical points.
\begin{Prop}
The critical points of \eqref{energy-functional} are given by
\begin{align*}
\delta_2\Delta\tau(\phi)=\delta_2 R^N(d\phi(e_\alpha),\tau(\phi))d\phi(e_\alpha)+\delta_1\tau(\phi),
\end{align*}
where \(\tau(\phi):=\tr_h\nabla d\phi\) is the tension field of the map \(\phi\).
\end{Prop}
\begin{proof}
We choose Riemannian normal coordinates that satisfy \(\nabla_{\partial_t}e_\alpha=0\) at the respective point.
Consider a variation of the map \(\phi\), that is \(\phi_t\colon (-\epsilon,\epsilon)\times M\to N\), which 
satisfies \(\frac{\nabla}{\partial t}\phi\big|_{t=0}=\eta\).
It is well-known that
\begin{align*}
\frac{d}{dt}\big|_{t=0}\int_M|d\phi_t|^2\dv=-2\int_M\langle\eta,\tau(\phi)\rangle\dv.
\end{align*}
In addition, we find
\begin{align*}
\frac{d}{dt}\big|_{t=0}\int_M|\tau(\phi_t)|^2\dv=&2\int_M\langle\frac{\nabla}{\partial t}\nabla_{e_\alpha}d\phi_t(e_\alpha),\tau(\phi_t)\rangle\dv\big|_{t=0} \\
=&2\int_M(\langle R^N(d\phi_t(\partial_t),d\phi_t(e_\alpha))d\phi_t(e_\alpha),\tau(\phi_t)\rangle \\
&+\langle\nabla_{e_\alpha}\nabla_{e_\alpha}d\phi_t(\partial_t),\tau(\phi_t)\rangle)\dv\big|_{t=0} \\
=&2\int_M\langle\eta,\Delta\tau(\phi)-R^N(d\phi(e_\alpha),\tau(\phi))d\phi(e_\alpha)\rangle\dv.
\end{align*}
Adding up both contributions yields the claim.
\end{proof}

Solutions of \eqref{euler-lagrange} will be called \emph{interpolating sesqui-harmonic maps}.
\footnote{Finding an appropriate name for solutions of \eqref{euler-lagrange} turned out to be subtle. The author would like to thank John Wood for suggesting the word ``sesqui''.}

\begin{Bem}
Choosing \(\delta_1=2,\delta_2=1\) and \(M=S^4,N=S^k\) solutions
of \eqref{euler-lagrange} were called \emph{quasi-biharmonic maps} in \cite{MR2094320}.
These arise when considering a sequence of weakly intrinsic biharmonic
maps in dimension four. When taking the limit, one finds that quasi-biharmonic spheres separate
at finitely many points as in many conformally invariant variational problems.
\end{Bem}

\begin{Bem}
There is another way how we can think of \eqref{euler-lagrange}.
If we interpret the biharmonic map equation as acting with the Jacobi-field 
operator \(J\) on the tension field \(\tau(\phi)\), then we may rewrite
the equation for interpolating sesqui-harmonic maps as
\begin{align*}
J(\tau(\phi))=\frac{\delta_1}{\delta_2}\tau(\phi).
\end{align*}
Since the Jacobi-field operator is elliptic it has a discrete spectrum whenever
\(M\) is a closed manifold. In this case the equation for interpolating sesqui-harmonic maps
can be thought of as an eigenvalue equation for the Jacobi-field operator.
\end{Bem}

\begin{Bem}
In order to highlight the dependence of the action functional \eqref{energy-functional} on
the metric on the domain \(M\) we write
\begin{align*}
E_{\delta_1,\delta_2}(\phi,h)=\delta_1\int_M|d\phi|_h^2 \dv_{h}+\delta_2\int_M|\tau_h(\phi)|_h^2\dv_{h},
\end{align*}
where \(\dv_h\) represents the volume element of the metric \(h\).
If we perform a rescaling of the metric by a constant factor \(\tilde h:=\lambda^2h\),
the action functional transforms as
\begin{align*}
E_{\delta_1,\delta_2}(\phi,\tilde h)=\delta_1\int_M|d\phi|_h^2\lambda^{n-2}\dv_h + \delta_2\int_M|\tau_h(\phi)|_h^2\lambda^{n-4}\dv_h.
\end{align*}
This clearly reflects the fact that the action functional for harmonic maps is scale-invariant in two dimensions whereas
the action functional for biharmonic maps is scale-invariant in four dimensions. 
We can conclude that the action for interpolating sesqui-harmonic maps is not scale-invariant in any dimension and
we may expect that interpolating sesqui-harmonic maps may be most interesting if \(\dim M=2,3,4\).
\end{Bem}

In order to highlight the analytical structure of \eqref{euler-lagrange}
we take a look at the case of a spherical target.
For biharmonic maps this was carried out in \cite{MR2375314} making use of a different method.

\begin{Prop}
For \(\phi\colon M\to S^n\subset\R^{n+1}\) with the constant curvature metric,
the equation for interpolating sesqui-harmonic maps \eqref{euler-lagrange} acquires the form 
\begin{align}
\delta_2(\Delta^2\phi+(|\Delta\phi|^2+\Delta|d\phi|^2+2\langle d\phi,\nabla\Delta\phi\rangle+2|d\phi|^4)\phi+2\nabla(|d\phi|^2d\phi))=\delta_1(\Delta\phi+|d\phi|^2\phi).
\end{align}

\end{Prop}
\begin{proof}
Recall that for a spherical target the tension field has the simple form
\begin{align*}
\tau(\phi)=\Delta\phi+|d\phi|^2\phi.
\end{align*}

Since we assume that \(N=S^n\) with constant curvature the term on the right hand side of \eqref{euler-lagrange}
acquires the form
\begin{align*}
-R^N(d\phi(e_\alpha),\tau(\phi))d\phi(e_\alpha)=&|d\phi|^2\tau(\phi)-\langle d\phi(e_\alpha),\tau(\phi)\rangle d\phi(e_\alpha)\\
=&|d\phi|^2\Delta\phi+|d\phi|^4\phi-\langle d\phi(e_\alpha),\Delta\phi\rangle d\phi(e_\alpha).
\end{align*}
Using the special structure of the Levi-Civita connection on \(S^n\subset \R^{n+1}\) we calculate
\begin{align*}
\Delta\tau(\phi)=&\nabla\nabla(\Delta\phi+|d\phi|^2\phi) \\
=&\nabla(\nabla\Delta\phi+\langle d\phi,\Delta\phi\rangle\phi+\nabla(|d\phi|^2\phi))\\
=&\Delta^2\phi+\langle d\phi,\nabla\Delta\phi\rangle\phi
+\nabla(\langle d\phi,\Delta\phi\rangle\phi)
+\Delta(|d\phi|^2\phi)+\langle d\phi,\nabla(|d\phi|^2\phi)\rangle\phi\\
=&\Delta^2\phi+|\Delta\phi|^2\phi+2\langle d\phi,\nabla\Delta\phi\rangle\phi
+\langle d\phi,\Delta\phi\rangle d\phi+\Delta(|d\phi|^2\phi)+|d\phi|^4\phi.
\end{align*}
Combining both equations yields the claim.
\end{proof}

By varying \eqref{energy-functional} with respect to the domain metric
we obtain the energy-momentum tensor. Since the energy-momentum tensor
for both harmonic and biharmonic maps is well-known in the literature
we can directly give the desired result.

\begin{Prop}
The energy-momentum tensor associated to \eqref{energy-functional} is given by
\begin{align}
\label{energy-momentum}
T(X,Y)=&\delta_1(\langle d\phi(X),d\phi(Y)\rangle-\frac{1}{2}|d\phi|^2h(X,Y)) \\
\nonumber&+\delta_2\big(\frac{1}{2}|\tau(\phi)|^2h(X,Y)+\langle d\phi,\nabla\tau(\phi)\rangle h(X,Y) \\
\nonumber&-\langle d\phi(X),\nabla_Y\tau(\phi)\rangle-\langle d\phi(Y),\nabla_X\tau(\phi)\rangle\big),
\end{align}
where \(X,Y\) are vector fields on \(M\).
\end{Prop}

\begin{proof}
We consider a variation of the metric on \(M\), that is 
\begin{align*}
\frac{d}{dt}\big|_{t=0}h_t=k,
\end{align*}
where \(k\) is a symmetric \((2,0)\)-tensor. The energy-momentum tensor for harmonic maps can be computed as \cite{MR655417}
\begin{align*}
\frac{d}{dt}\big|_{t=0}\int_M|d\phi|^2\dv_{h_t}
=\int_M k^{\alpha\beta}\big(\langle d\phi(e_\alpha),d\phi(e_\beta)\rangle-\frac{1}{2}|d\phi|^2h_{\alpha\beta}\big)\dv_h.
\end{align*}
Deriving the energy-momentum tensor for biharmonic maps is more involved as one has to vary the connection
of the domain since it depends on the metric. The energy-momentum tensor for biharmonic maps was already presented in \cite{MR886529},
a rigorous derivation was obtained in \cite[Theorem 2.4]{MR2395125}, that is
\begin{align*}
\frac{d}{dt}\big|_{t=0}\int_M|\tau(\phi)|^2\dv_{h_t}
=-\int_M &k^{\alpha\beta}(\frac{1}{2}|\tau(\phi)|^2h_{\alpha\beta}+\langle d\phi,\nabla\tau(\phi)\rangle h_{\alpha\beta} \\
&-\langle d\phi(e_\alpha),\nabla_{e_\beta}\tau(\phi)\rangle-\langle d\phi(e_\beta),\nabla_{e_\alpha}\tau(\phi)\rangle
) \dv_h.
\end{align*}
Combining both formulas concludes the proof.
\end{proof}

It can be directly seen that the energy-momentum tensor \eqref{energy-momentum} is symmetric.
For the sake of completeness we prove the following
\begin{Prop}
The energy-momentum tensor \eqref{energy-momentum} is divergence-free.
\end{Prop}
\begin{proof}
We choose a local orthonormal basis \(e_\alpha\) and set
\begin{align*}
T_{\alpha\beta}:=T(e_\alpha,e_\beta)=&\delta_1(\langle d\phi(e_\alpha),d\phi(e_\beta)\rangle-\frac{1}{2}|d\phi|^2h_{\alpha\beta}) \\
\nonumber&+\delta_2\big(\frac{1}{2}|\tau(\phi)|^2h_{\alpha\beta}+\langle d\phi,\nabla\tau(\phi)\rangle h_{\alpha\beta} \\
\nonumber&-\langle d\phi(e_\alpha),\nabla_{e_\beta}\tau(\phi)\rangle-\langle d\phi(e_\beta),\nabla_{e_\alpha}\tau(\phi)\rangle\big).
\end{align*}
By a direct calculation we find
\begin{align*}
\nabla_{e_\alpha}(\langle d\phi(e_\alpha),d\phi(e_\beta)\rangle-\frac{1}{2}|d\phi|^2h_{\alpha\beta})=\langle\tau(\phi),d\phi(e_\beta)\rangle
\end{align*}
and also
\begin{align*}
\nabla_{e_\alpha}&(\frac{1}{2}|\tau(\phi)|^2h_{\alpha\beta}+\langle d\phi,\nabla\tau(\phi)\rangle h_{\alpha\beta}
-\langle d\phi(e_\alpha),\nabla_{e_\beta}\tau(\phi)\rangle-\langle d\phi(e_\beta),\nabla_{e_\alpha}\tau(\phi)\rangle) \\
=&\langle\nabla_{e_\beta}\tau(\phi),\tau(\phi)\rangle+\langle\nabla_{e_\beta}d\phi(e_\gamma),\nabla_{e_\gamma}\tau(\phi)\rangle
+\langle d\phi(e_\gamma),\nabla_{e_\beta}\nabla_{e_\gamma}\tau(\phi)\rangle \\
&-\langle\tau(\phi),\nabla_{e_\beta}\tau(\phi)\rangle
-\langle d\phi(e_\alpha),\nabla_{e_\alpha}\nabla_{e_\beta}\tau(\phi)\rangle
-\langle\nabla_{e_\alpha}d\phi(e_\beta),\nabla_{e_\alpha}\tau(\phi)\rangle
-\langle d\phi(e_\beta),\nabla\nabla\tau(\phi)\rangle \\
=&\langle d\phi(e_\alpha),R^N(d\phi(e_\beta),d\phi(e_\alpha))\tau(\phi)\rangle-\langle d\phi(e_\beta),\Delta\tau(\phi)\rangle,
\end{align*}
where we used the torsion-freeness of the Levi-Civita connection.
Adding up both equations yields the claim.
\end{proof}

\subsection{Conservation laws for targets with symmetries}
In this subsection we discuss how to obtain a conservation law for solutions of
the interpolating sesqui-harmonic map equation in the case that the target manifold has a certain
amount of symmetry, more precisely if it possesses Killing vector fields.
A similar discussion has been performed in \cite{MR3735550}.

To this end let \(\xi\) be a diffeomorphism that generates a one-parameter family of vector fields \(X\).
Then we know that
\begin{align*}
\frac{d}{dt}\big|_{t=0}\xi^\ast g=\cL_Xg,
\end{align*}
where \(\cL\) denotes the Lie-derivative acting on the metric.
In terms of local coordinates the Lie-derivative of the metric is given by
\begin{align*}
\cL_Xg_{ij}=\nabla_iX_j+\nabla_jX_i.
\end{align*}

This enables us to give the following 
\begin{Dfn}
Let \(\xi\) be a diffeomorphism that generates a one-parameter family of vector fields \(X\) on \(N\).
We say that \(X\) generates a symmetry for the action \(E_{\delta_1,\delta_2}(\phi,\xi^\ast g)\) if
\begin{align*}
\frac{d}{dt}\big|_{t=0}E_{\delta_1,\delta_2}(\phi,\xi^\ast g)=\int_M\cL_X(\delta_1|d\phi|^2+\delta_2|\tau(\phi)|^2)\dv=0,
\end{align*}
where the Lie-derivative is acting on the metric \(g\).
\end{Dfn}

Note that if \(X\) generates an isometry then \(\cL_Xg=0\) such that we have to require the existence of Killing
vector fields on the target.

In the following we will make use of the following facts:
\begin{Lem}
If \(X\) is a Killing vector field on the target \(N\), then  
\begin{align}
\label{lie-curvature}\nabla^2_{Y,Z}X&=-R^N(X,Y)Z,\\
\label{lie-christoffel}\cL_X\Gamma^k_{ij}&=\nabla_i\nabla_jX^k-R^k_{~ijl}X^l,
\end{align}
where \(\Gamma^k_{ij}\) are the Christoffel symbols on \(N\), \(R^k_{~ijl}\) the components of 
the Riemannian curvature tensor on \(N\) and \(Y,Z\) vector fields on \(N\).
\end{Lem}

\begin{Lem}
Let \(\phi\colon M\to N\) be a smooth solution of \eqref{euler-lagrange} and assume that \(N\)
admits a Killing vector field \(X\).
Then the Lie-derivative acting on the metric \(g\) of the energy density is given by
\begin{align*}
\cL_X(\delta_1|d\phi|^2+\delta_2|\tau(\phi)|^2)
=&2\delta_1\nabla_{e_\alpha}\langle d\phi(e_\alpha),X(\phi)\rangle \\
\nonumber&+2\delta_2\nabla_{e_\alpha}(\langle\tau(\phi),\nabla_{e_\alpha}X(\phi)\rangle-\delta_2\langle\nabla_{e_\alpha}\tau(\phi),X(\phi)\rangle).
\end{align*}
\end{Lem}
\begin{proof}
We choose Riemannian normal coordinates \(x_\alpha\) on \(M\) and calculate
\begin{align*}
\cL_X|d\phi|^2=&(\cL_Xg)_{ij}\frac{\partial\phi^i}{\partial x^\alpha}\frac{\partial\phi^j}{\partial x^\beta}h^{\alpha\beta} \\
=&2\nabla_iX_j\frac{\partial\phi^i}{\partial x^\alpha}\frac{\partial\phi^j}{\partial x^\beta}h^{\alpha\beta} \\
=&2\langle d\phi(e_\alpha),\nabla_{e_\alpha}(X(\phi))\rangle \\
=&2\nabla_{e_\alpha}\langle d\phi(e_\alpha),X(\phi)\rangle-2\langle\tau(\phi),X(\phi)\rangle.
\end{align*}
To calculate the variation of the tension field with respect to the target metric
we first of all note that
\begin{align*}
\cL_X|\tau(\phi)|^2=&2\tau^i(\phi)\tau^j(\phi)\nabla_iX_j+2g_{ij}\tau^i(\phi)\cL_X\tau^j(\phi).
\end{align*}

Making use of the local expression of the tension field this yields
\begin{align*}
\cL_X\tau^j(\phi)=&\cL_X(\Delta\phi^j+h^{\alpha\beta}\frac{\partial\phi^k}{\partial x^\alpha}\frac{\partial\phi^l}{\partial x^\beta}\Gamma^j_{kl})\\
=&h^{\alpha\beta}\frac{\partial\phi^k}{\partial x^\alpha}\frac{\partial\phi^l}{\partial x^\beta}\cL_X\Gamma^j_{kl} \\
=&h^{\alpha\beta}\frac{\partial\phi^k}{\partial x^\alpha}\frac{\partial\phi^l}{\partial x^\beta}(\nabla_k\nabla_lX^j-R^j_{~klr}X^r),
\end{align*}
where we used \eqref{lie-christoffel} in the last step.

This allows us to infer
\begin{align*}
g_{ij}\tau^i(\phi)\cL_X\tau^j(\phi)=\langle\tau(\phi),\nabla_{d\phi(e_\alpha)}\nabla_{e_\alpha}X(\phi)\rangle+\langle R^N(\tau(\phi),d\phi(e_\alpha))d\phi(e_\alpha),X\rangle.
\end{align*}

Combining both equations we find 
\begin{align*}
\cL_X|\tau(\phi)|^2=&2\nabla_{e_\alpha}\langle\tau(\phi),\nabla_{e_\alpha}X(\phi)\rangle
-2\nabla_{e_\alpha}\langle\nabla_{e_\alpha}\tau(\phi),X(\phi)\rangle
+2\langle\Delta\tau(\phi),X(\phi)\rangle \\
&+2\langle R^N(\tau(\phi),d\phi(e_\alpha))d\phi(e_\alpha),X\rangle.
\end{align*}
The result follows by adding up both contributions.
\end{proof}

\begin{Prop}
Let \(\phi\colon M\to N\) be a smooth solution of \eqref{euler-lagrange} and assume that \(N\)
admits a Killing vector field \(X\).
Then the following vector field is divergence free
\begin{align}
\label{noether}
J_\alpha:=\delta_1\langle d\phi(e_\alpha),X(\phi)\rangle
+\delta_2\langle\tau(\phi),\nabla_{e_\alpha}X(\phi)\rangle-\delta_2\langle\nabla_{e_\alpha}\tau(\phi),X(\phi)\rangle.
\end{align}

\end{Prop}
\begin{proof}
A direct calculation yields
\begin{align*}
\nabla_{e_\alpha}J_\alpha=&\delta_1\langle\tau(\phi),X(\phi)\rangle
+\delta_1\underbrace{\langle d\phi(e_\alpha),\nabla_{e_\alpha} X(\phi)\rangle}_{=0}
+\delta_2\langle\tau(\phi),\nabla^2_{e_\alpha,e_\alpha} X(\phi)\rangle-\delta_2\langle\Delta\tau(\phi),X(\phi)\rangle\\
=&\langle X(\phi),\delta_1\tau(\phi)+\delta_2R^N(d\phi(e_\alpha),\tau(\phi))d\phi(e_\alpha)-\delta_2\Delta\tau(\phi)\rangle\\
=&0,
\end{align*}
where we used \eqref{lie-curvature} and that \(\phi\) is a solution of \eqref{euler-lagrange} in the last step.
\end{proof}

\begin{Bem}
In the physics literature the vector field \eqref{noether} is usually called 
\emph{Noether current}.
\end{Bem}

\section{Explicit solutions of the interpolating sesqui-harmonic map equation}
In this section we want to derive several explicit solutions to the 
Euler-Lagrange equation \eqref{euler-lagrange}.
We can confirm that solutions may have a different behavior than biharmonic or
harmonic maps.

Let us start in the most simple setup possible.
\begin{Bsp}
Suppose that \(M=N=S^1\) and by \(s\) we denote the global coordinate on \(S^1\).
Then \eqref{euler-lagrange} acquires the form
\begin{align*}
\delta_2\phi^{''''}(s)=\delta_1\phi^{''}(s).
\end{align*}
Taking an ansatz of the form 
\begin{align*}
\phi(s)=\sum_k a_ke^{iks}
\end{align*}
we obtain 
\begin{align*}
\sum_k a_kk^2(\delta_2k^2-\delta_1)e^{iks}=0.
\end{align*}
Consequently, we have to impose the condition \(k^2=\frac{\delta_1}{\delta_2}\).
In particular, this shows that there does not exist a solution of \eqref{euler-lagrange}
on \(S^1\) if \(\delta_1\) and \(\delta_2\) have opposite sign.
\end{Bsp}

\begin{Bsp}
Consider the case \(M=\R^2\) and \(N=\R\). In this case being interpolating sesqui-harmonic
means to find a function \(f\colon M\to N\) that solves 
\begin{align*}
\delta_2(\partial^2_x+\partial^2_y)^2f=\delta_1(\partial^2_x+\partial^2_y)f,
\end{align*}
where \(x,y\) denote the canonical coordinates in \(\R^2\).
If we make a separation ansatz of the form 
\[f(x,y)=e^{\alpha x}e^{\beta y}\]
then we are lead to the following algebraic expression
\begin{align*}
\delta_2(\alpha^2+\beta^2)=\delta_1.
\end{align*}
Let us distinguish the following cases
\begin{enumerate}
 \item \(\delta_1=0\), that is \(f\) is biharmonic. In this case \(\alpha^2+\beta^2=0\) and
  we have to choose \(\alpha,\beta\in\C\).
 \item \(\delta_2=0\), that is \(f\) is harmonic. In this case there are no restrictions on \(\alpha,\beta\).
 \item If \(\delta_1\delta_2>0\) we have to impose the condition \(\alpha^2+\beta^2>0\) meaning that \(\alpha,\beta\in\R\).
 \item If \(\delta_1\delta_2<0\) we find that \(\alpha^2+\beta^2<0\) meaning that \(\alpha,\beta\in\C\).
\end{enumerate}
This again shows that interpolating sesqui-harmonic functions may be very different from both harmonic and biharmonic functions.
\end{Bsp}

\subsection{Interpolating sesqui-harmonic functions in flat space}
In this section we study interpolating sesqui-harmonic functions in flat space.

First, suppose that \(M=N=\R\) and we denote the global coordinate on \(\R\) by \(x\).
Then \eqref{euler-lagrange} acquires the form
\begin{align*}
\delta_2\phi^{''''}(x)=\delta_1\phi^{''}(x).
\end{align*}
This can be integrated as
\begin{align*}
\phi(x)=\frac{\delta_2}{\delta_1}(c_1e^{\sqrt{\frac{\delta_1}{\delta_2}}x}+c_2e^{-\sqrt{\frac{\delta_1}{\delta_2}}x})+c_3x+c_4,
\end{align*}
where \(c_i,i=1,\ldots 4\) are integration constants.
It is interesting to note that both limits \(\delta_1\to 0\) and \(\delta_2\to 0\) do not exist.
Consequently, the solution of the interpolating sesqui-harmonic function equation does neither reduce
to a solution of the harmonic or the biharmonic function equation.
In the following we will analyze if the same behavior persists in higher dimensions.

If we take \(M=(\R^n,\delta)\) and \(N=(\R,\delta)\)
both with the Euclidean metric the equation for
interpolating sesqui-harmonic maps turns into
\begin{align}
\label{interpolating sesqui-harmonic-rn}
\delta_2\Delta\Delta f=\delta_1\Delta f,
\end{align}
where \(f\colon\R^n\to\R\).
Although this equation is linear we may expect some analytical difficulties
since we do not have a maximum principle available for fourth order equations.

A full in detail analysis of this equation is far beyond the scope of this article.
Nevertheless, we will again see that solutions of \eqref{interpolating sesqui-harmonic-rn} may 
be very different from harmonic and biharmonic functions.
We will be looking for radial solutions of \eqref{interpolating sesqui-harmonic-rn}, where \(r:=\sqrt{x_1^2+\ldots+x_n^2}\).
In this case the Laplacian has the form
\begin{align*}
\Delta=\frac{d^2}{dr^2}+\frac{n-1}{r}\frac{d}{dr}.
\end{align*}
Recall that the fundamental solution of the Laplace equation in \(\R^n\) is given by
\begin{align*}
H_{\Delta}(x,y)=
\begin{cases}
|x-y|^{2-n},& n\geq 3, \\ \log|x-y|,& n=2,
\end{cases}
\end{align*}
whereas for biharmonic functions it acquires the form
\begin{align*}
H_{\Delta^2}(x,y)=
\begin{cases}
|x-y|^{4-n},& n\geq 5, \\ \log|x-y|,& n=4,\\|x-y|,& n=3.
\end{cases}
\end{align*}
Note that we did not write down any normalization of the fundamental solutions.

We cannot expect to find a unique solution to \eqref{interpolating sesqui-harmonic-rn} since
we can always add a harmonic function once we have constructed a solution to \eqref{interpolating sesqui-harmonic-rn}.
Since we are considering \(\R^n\) instead of a curved manifold at the moment all curvature terms in \eqref{euler-lagrange}
vanish and we are dealing with a linear problem.

Instead of trying to directly solve \eqref{interpolating sesqui-harmonic-rn} we rewrite the equation as follows
\begin{align*}
\Delta(\Delta f-\frac{\delta_1}{\delta_2}f)=0.
\end{align*}
Assume that \(n\geq 3\) and making use of the fundamental solution of the Laplacian
we may solve
\begin{align*}
\Delta f(r)-\frac{\delta_1}{\delta_2}f(r)=r^{2-n},
\end{align*}
which then provides an interpolating sesqui-harmonic function. This yields the following 
ordinary differential equation
\begin{align}
\label{flatspace-ode}
f''(r)+\frac{n-1}{r}f'(r)-\frac{\delta_1}{\delta_2}f(r)=r^{2-n}.
\end{align}
This equation can be solved explicitly in terms of a linear combination of Bessel functions in any dimension.
Since the general solution is rather lengthy, we only give some explicit solutions for a fixed dimension.
\begin{itemize}
 \item Suppose that \(n=3\) then the solution of \eqref{flatspace-ode} is given by
 \begin{align*}
  f(r)=c_1\frac{e^{-\sqrt{\frac{\delta_1}{\delta_2}}r}}{r}+c_2\frac{e^{\sqrt{\frac{\delta_1}{\delta_2}}r}}{\sqrt{\delta_1/\delta_2}r}-\frac{\delta_2}{\delta_1}\frac{1}{r}.
 \end{align*}
 As in the one-dimensional case both limits \(\delta_1\to 0\) and \(\delta_2\to 0\) do not exist.
\item Suppose that \(n=4\) and that \(\frac{\delta_1}{\delta_2}>0\). Then the solution of \eqref{flatspace-ode} is given by
\begin{align*}
f(r)=&c_1\frac{J_1(\sqrt{\frac{\delta_1}{\delta_2}}r)}{r}+c_2\frac{Y_1(\sqrt{\frac{\delta_1}{\delta_2}}r)}{r} \\
&+\frac{\pi}{2\sqrt{\delta_1/\delta_2}r}(J_1(\sqrt{\frac{\delta_1}{\delta_2}}r)Y_0(\sqrt{\frac{\delta_1}{\delta_2}}r)
-J_0(\sqrt{\frac{\delta_1}{\delta_2}}r)Y_1(\sqrt{\frac{\delta_1}{\delta_2}}r)).
\end{align*}
If \(\frac{\delta_1}{\delta_2}<0\) then we obtain the solution by an analytic continuation.
\item The qualitative behavior of solutions to \eqref{flatspace-ode} for \(n\geq 5\) seems to be the same as above.
\end{itemize}
It becomes obvious that solutions of \eqref{interpolating sesqui-harmonic-rn} may show
a different qualitative behavior compared to biharmonic functions.
Moreover, as one should expect, the qualitative behavior depends heavily
on the sign of the product \(\delta_1\delta_2\).

\subsection{Interpolating sesqui-harmonic curves on the three-dimensional sphere}
In this subsection we study interpolating sesqui-harmonic curves on three-dimensional spheres
with the round metric, where we follow the ideas from \cite{MR1863283}.

To this end let \((N,g)\) be a three-dimensional Riemannian manifold with constant sectional curvature \(K\).
Moreover, let \(\gamma\colon I\to N\) be a smooth curve that is parametrized by arc length.
Let \(\{T,N,B\}\) be an orthonormal frame field of \(TN\) along the curve \(\gamma\).
Here, \(T=\gamma'\) is the unit tangent vector of \(\gamma\), \(N\) the unit normal field
and \(B\) is chosen such that \(\{T,N,B\}\) forms a positive oriented basis.

In this setup we have the following Frenet equations for the curve \(\gamma\)
\begin{align}
\label{frenet-three-dimensions}
\nabla_TT=k_gN,\qquad \nabla_TN=-k_gT+\tau_gB,\qquad \nabla_TB=-\tau_gN.
\end{align}

\begin{Lem}
Let \(\gamma\colon I\to N\) be a curve in a three-dimensional Riemannian manifold.
Then the curve \(\gamma\) is interpolating sesqui-harmonic if the following equation holds
\begin{align*}
(-3\delta_2k_gk_g')T+(\delta_2(k_g''-k_g^3-k_g\tau_g^2+k_gK)-\delta_1k_g)N
+\delta_2(2k_g'\tau_g+k_g\tau'_g)B=0.
\end{align*}
\end{Lem}
\begin{proof}
Making use of the Frenet equations \eqref{frenet-three-dimensions} a direct calculation yields
\begin{align*}
\nabla^3_TT=(-3k_gk_g')T+(k_g''-k_g^3-k_g\tau^2_g)N+(2k_g'\tau_g+k_g\tau_g')B=0.
\end{align*}
Using that the sectional curvature of \(N\) is given by \(K=R^N(T,N,N,T)\) we obtain the claim.
\end{proof}

\begin{Cor}
Let \(\gamma\colon I\to N\) be a curve in a three-dimensional Riemannian manifold.
Then the curve \(\gamma\) is interpolating sesqui-harmonic if the following system holds
\begin{align*}
k_gk_g'=0,\qquad 2k_g'\tau_g+k_g\tau'_g=0,\qquad \delta_2(k_g''-k_g^3-k_g\tau_g^2+k_gK)=\delta_1k_g.
\end{align*}
The non-geodesic solutions \((k_g\neq 0)\) of this system are given by
\begin{align}
\label{curve-three-semibiharmonic}
k_g=const\neq 0,\qquad \tau_g=const,\qquad \delta_2(k_g^2+\tau_g^2)=\delta_2K-\delta_1.
\end{align}
\end{Cor}

We directly obtain the following characterization of interpolating sesqui-harmonic curves:
\begin{Prop}
\begin{enumerate}
 \item Let \(\gamma\colon I\to N\) be a curve in a three-dimensional Riemannian manifold.
If \(K\leq\frac{\delta_1}{\delta_2}\) then any interpolating sesqui-harmonic curve is a geodesic.
\item To obtain a non-geodesic interpolating sesqui-harmonic curve \(\gamma\colon I\to S^3\) we have to demand that 
\(\delta_2>\delta_1\).
\end{enumerate}
\end{Prop}

\begin{Prop}
Let \(\gamma\colon I\to S^3\) be a curve on the three-dimensional sphere with the round metric.
The curve \(\gamma\) is interpolating sesqui-harmonic if the following equation holds
\begin{align}
\label{curve-s3}
\gamma^{''''}+(1-\delta_1+\delta_2)\gamma^{''}+(-k_g^2-\delta_1+\delta_2)\gamma=0.
\end{align}
\end{Prop}
\begin{proof}
Differentiating the first equation of \eqref{frenet-three-dimensions} we find
\begin{align*}
\nabla^2_TN=&-k_g'T-k_g\nabla_TT+\tau'_gB+\tau_g\nabla_TB \\
=&-k_g\nabla_TT+\tau_g\nabla_TB \\
=&-(k_g^2+\tau_g^2)N \\
=&(-\delta_2+\delta_1)N,
\end{align*}
where we used \eqref{curve-three-semibiharmonic}.
Moreover, employing the formula for the Levi-Civita connection on \(S^3\subset\R^4\)
\begin{align*}
\nabla_TX=X'+\langle T,X\rangle\gamma
\end{align*}
we get the equations
\begin{align*}
\nabla^2_TN=&\nabla_T(N'+\langle T,N\rangle\gamma) 
=N''+\langle T,\nabla_TN\rangle\gamma 
=N''-k_g\gamma, \\
\nabla_TT=&k_gN=\gamma''+|\gamma'|^2\gamma=\gamma''+\gamma.
\end{align*}
Combining the equations for \(\nabla^2_TN\) we obtain
\begin{align*}
N''-k_g\gamma=(\delta_1-\delta_2)N
\end{align*}
and rewriting this as an equation for \(\gamma\) yields the claim.
\end{proof}

\begin{Prop}
Let \(\gamma\colon I\to S^3\) be a curve on the three-dimensional sphere with the round metric.
If \(k_g^2=\delta_2-\delta_1\) the interpolating sesqui-harmonic curves are given by
\begin{align}
\label{solution-s3-a}
\gamma(t)=\big(\frac{\cos(\sqrt{(1-\delta_1+\delta_2)}t)}{\sqrt{(1-\delta_1+\delta_2)}},\frac{\sin(\sqrt{(1-\delta_1+\delta_2)}t)}{\sqrt{(1-\delta_1+\delta_2)}},d_1,d_2\big)
\end{align}
where \(\frac{1}{1-\delta_1+\delta_2}+d_1^2+d_2^2=1\).
\end{Prop}

\begin{proof}
Making use of the assumptions \eqref{curve-s3} simplifies as
\begin{align*}
\gamma^{''''}+(1-\delta_1+\delta_2)\gamma^{''}=0.
\end{align*}
Solving this differential equation together with 
the constraints \(|\gamma|^2=1\) and \(|\gamma'|^2=1\) yields the claim.
\end{proof}

\begin{Bem}
Note that it is required in \eqref{solution-s3-a} that
\(1-\delta_1+\delta_2>1\), which is equivalent to \(\delta_2>\delta_1\).
This is consistent with the assumption \(k_g^2=\delta_2-\delta_1\).
\end{Bem}

\begin{Satz}
Let \(\gamma\colon I\to S^3\) be a curve on the three-dimensional sphere with the round metric
and suppose that \(\delta_2>\delta_1\).
Then the non-geodesic solution to \eqref{curve-s3} is given by
\begin{align}
\label{curve-s3-general-solution}
\gamma(t)=\frac{1}{\sqrt{a_1^2-a_2^2}}\big(\sqrt{1-a_2^2}\cos(a_1t),\sqrt{1-a_2^2}\sin(a_1t),\sqrt{a_1^2-1}\cos(a_2t),\sqrt{a_1^2-1}\sin(a_2t)\big)
\end{align}
with the constants 
\begin{align*}
a_1&:=\frac{1}{\sqrt{2}}\sqrt{1-\delta_1+\delta_2+\sqrt{(1+\delta_1-\delta_2)^2+4k_g^2}}, \\
a_2&:=\frac{1}{\sqrt{2}}\sqrt{1-\delta_1+\delta_2-\sqrt{(1+\delta_1-\delta_2)^2+4k_g^2}}.
\end{align*}
\end{Satz}

\begin{proof}
The most general ansatz for a solution of \eqref{curve-s3} is given by
\begin{align*}
\gamma(t)=c_1\cos(at)+c_2\sin(at)+c_3\cos(bt)+c_4\sin(bt),
\end{align*}
where \(c_i,i=1\ldots 4\) are mutually perpendicular and \(a\) and \(b\) are real numbers.
This leads to the following quadratic equation for both \(a\) and \(b\)
\begin{align*}
a^4-a^2(1-\delta_1+\delta_2)+(-k_g^2-\delta_1+\delta_2)=0.
\end{align*}
We obtain the two solutions
\begin{align*}
a_1^2&=\frac{1}{2}\big(1-\delta_1+\delta_2+\sqrt{(1+\delta_1-\delta_2)^2+4k_g^2}\big), \\
a_2^2&=\frac{1}{2}\big(1-\delta_1+\delta_2-\sqrt{(1+\delta_1-\delta_2)^2+4k_g^2}\big).
\end{align*}
Moreover, the constraints \(|\gamma|^2=1\) and \(|\gamma'|^2=1\) give the two equations
\begin{align*}
|c_1|^2+|c_2|^2=1,\qquad a_1^2|c_1|^2+a_2^2|c_2|^2=1.
\end{align*}
Solving this system for \(|c_1|^2\) and \(|c_2|^2\) yields the claim.
\end{proof}

\begin{Bem}
If we analyze the constants appearing in \eqref{curve-s3-general-solution} then we find
that we have to demand the condition \(\delta_2-\delta_1>k_g^2>0\) in order 
to obtain a real-valued constant \(a_2\).
In addition, we find
\begin{align*}
a_1^2-a_2^2&=\sqrt{(1+\delta_1-\delta_2)^2+4k_g^2}>0, \\
1-a_2^2&=\frac{1}{2}(1+\delta_1-\delta_2+\sqrt{(1+\delta_1-\delta_2)^2+4k_g^2})>0,\\
a_1^2-1&=\frac{1}{2}(-1-\delta_1+\delta_2+\sqrt{(1+\delta_1-\delta_2)^2+4k_g^2})>0
\end{align*}
such that we do not get any further restrictions.
We conclude that we get a solution of \eqref{curve-three-semibiharmonic}
for all \(\delta_1,\delta_2\) satisfying \(\delta_2-\delta_1>k_g^2>0\).
\end{Bem}

\begin{Bem}
If we compare our results with \cite[Theorem 3.3]{MR1863283}, then we find that interpolating sesqui-harmonic curves on \(S^3\)
have the same qualitative behavior as biharmonic curves.
More precisely, we have the following two cases
\begin{enumerate}
 \item If \(k_g^2=\delta_2-\delta_1\), then \(\gamma\) is a circle of radius \(\frac{1}{\sqrt{1+k_g^2}}\).
 \item If \(\delta_2-\delta_1>k_g^2>0\), then \(\gamma\) is a geodesic 
  of the rescaled Clifford torus 
  \[
   S^1\big(\frac{\sqrt{(1+\delta_1-\delta_2+\sqrt{(1+\delta_1-\delta_2)^2+4k_g^2})}}{\sqrt{2}((1+\delta_1-\delta_2)^2+4k_g^2)^\frac{1}{4}}\big)
  \times S^1\big(\frac{\sqrt{-1-\delta_1+\delta_2+\sqrt{(1+\delta_1-\delta_2)^2+4k_g^2}}}{\sqrt{2}((1+\delta_1-\delta_2)^2+4k_g^2)^\frac{1}{4}}\big).
  \]
\end{enumerate}
Note that the solutions from above reduce to biharmonic curves (see \cite[Theorem 3.3]{MR1863283}) in the case of \(\delta_1=0,\delta_2=1\).

\end{Bem}

We want to close this subsection by mentioning that it is possible to generalize the results
obtained from above to higher-dimensional spheres as was done for biharmonic curves in \cite{MR1919374}.

\section{The qualitative behavior of solutions}
In this section we study the qualitative behavior of interpolating sesqui-harmonic maps.

In the case of a one-dimensional domain and the target being a Riemannian manifold,
the Euler-Lagrange equation reduces to
\begin{align}
\label{euler-lagrange-curve}
\delta_2\nabla^3_{\gamma'}\gamma'=\delta_2R^N(\gamma',\nabla_{\gamma'}\gamma')\gamma'+\delta_1\nabla_{\gamma'}\gamma',
\end{align}
where \(\gamma\colon I\to N\) and \(\gamma'\) denotes the derivative with respect to the curve parameter \(s\).
\begin{Prop}
Suppose that \(\gamma\colon I\to N\) is a smooth solution of \eqref{euler-lagrange}.
Then the following conservation type law holds
\begin{align*}
\big(\delta_2\frac{d^3}{ds^3}-\delta_1\frac{d}{ds}\big)\frac{1}{2}|\gamma'|^2=
\delta_2\frac{d}{ds}\frac{3}{2}|\nabla_{\gamma'}\gamma'|^2.
\end{align*}
\end{Prop}

\begin{proof}
We test \eqref{euler-lagrange-curve} with \(\gamma'\) and obtain
\begin{align*}
\delta_2\langle\nabla^3_{\gamma'}\gamma',\gamma'\rangle
=\delta_1\langle\nabla_{\gamma'}\gamma',\gamma'\rangle
=\frac{1}{2}\delta_1\frac{d}{ds}|\gamma'|^2.
\end{align*}
The left-hand side can be further simplified as
\begin{align*}
\langle\nabla^3_{\gamma'}\gamma',\gamma'\rangle=&\frac{d}{ds}\langle\nabla^2_{\gamma'}\gamma',\gamma'\rangle-\langle\nabla^2_{\gamma'}\gamma',\nabla_{\gamma'}\gamma'\rangle \\
=&\frac{d^2}{ds^2}\langle\nabla_{\gamma'}\gamma',\gamma'\rangle-\frac{3}{2}\frac{d}{ds}|\nabla_{\gamma'}\gamma'|^2 \\
=&\frac{d^3}{ds^3}\frac{1}{2}\langle\gamma',\gamma'\rangle-\frac{3}{2}\frac{d}{ds}|\nabla_{\gamma'}\gamma'|^2,
\end{align*}
which completes the proof.
\end{proof}

As already stated in the introduction it is obvious that harmonic maps solve \eqref{euler-lagrange}.
We will give several conditions under which interpolating sesqui-harmonic maps must be harmonic
generalizing several results from \cite{MR2124627,MR2004799}. 
To achieve these results we will frequently make use of the following Bochner formula:

\begin{Lem}
Let \(\phi\colon M\to N\) be a smooth solution of \eqref{euler-lagrange}.
Then the following Bochner formula holds
\begin{align}
\label{bochner-tension}
\Delta\frac{1}{2}|\tau(\phi)|^2=|\nabla\tau(\phi)|^2+\langle R^N(d\phi(e_\alpha),\tau(\phi))d\phi(e_\alpha),\tau(\phi)\rangle+\frac{\delta_1}{\delta_2}|\tau(\phi)|^2.
\end{align}
\end{Lem}
\begin{proof}
This follows by a direct calculation.
\end{proof}

\begin{Prop}
Suppose that \((M,h)\) is a compact Riemannian manifold.
Let \(\phi\colon M\to N\) be a smooth solution of \eqref{euler-lagrange}.
\begin{enumerate}
 \item If \(N\) has non-positive curvature \(K^N\leq 0\) and \(\delta_1,\delta_2\) have the same sign then \(\phi\) is harmonic.
 \item If \(|d\phi|^2\leq\frac{\delta_1}{|R^N|_{L^\infty}\delta_2}\) and \(\delta_1,\delta_2\) have the same sign then \(\phi\) is harmonic.
\end{enumerate}
\end{Prop}
\begin{proof}
The first statement follows directly from \eqref{bochner-tension} by application of the maximum principle.
For the second statement we estimate \eqref{bochner-tension} as
\begin{align*}
\Delta\frac{1}{2}|\tau(\phi)|^2=|\nabla\tau(\phi)|^2+(\frac{\delta_1}{\delta_2}-|R^N|_{L^\infty}|d\phi|^2)|\tau(\phi)|^2\geq 0
\end{align*}
due to the assumptions. The claim follows again due to the maximum principle.
\end{proof}

If we do not require \(M\) to be compact we can give the following result.
\begin{Prop}
Let \(\phi\colon M\to N\) be a Riemannian immersion that solves \eqref{euler-lagrange} with \(|\tau(\phi)|=const\).
If \(N\) has non-positive curvature \(K^N\leq 0\) and \(\delta_1,\delta_2\) have the same sign then \(\phi\) must be harmonic.
\end{Prop}

\begin{proof}
Via the maximum principle we obtain \(\nabla\tau(\phi)=0\) from \eqref{bochner-tension}.
By assumption the map \(\phi\) is an immersion such that
\begin{align*}
-|\tau(\phi)|^2=\langle d\phi,\nabla\tau(\phi)\rangle
\end{align*}
concluding the proof.
\end{proof}

In the case that \(\dim M=\dim N-1\) the assumption of \(N\) having
negative sectional curvature can be replaced by demanding negative Ricci curvature.

\begin{Satz}
Let \(\phi\colon M\to N\) be a Riemannian immersion.
Suppose that \(M\) is compact and \(\dim M=\dim N-1\).
If \(N\) has non-positive Ricci curvature and \(\delta_1,\delta_2\) have the same sign then \(\phi\) 
is interpolating sesqui-harmonic if and only if it is harmonic.
\end{Satz}
\begin{proof}
Since \(\phi\) is an immersion and \(\dim M=\dim N-1\) we obtain
\begin{align*}
R^N(d\phi(e_\alpha),\tau(\phi))d\phi(e_\alpha)=-\operatorname{Ric}^N(\tau(\phi)).
\end{align*}
Making use of \eqref{bochner-tension} we get
\begin{align*}
\Delta\frac{1}{2}|\tau(\phi)|^2=|\nabla\tau(\phi)|^2-\langle \operatorname{Ric}^N(\tau(\phi)),\tau(\phi)\rangle+\frac{\delta_1}{\delta_2}|\tau(\phi)|^2\geq 0
\end{align*}
due to the assumptions. The result follows by the maximum principle.
\end{proof}

As for harmonic maps (\cite[Theorem 2]{MR510549}) we can prove a unique continuation theorem for interpolating sesqui-harmonic maps.
To obtain this result we recall the following (\cite[p.248]{MR0092067})
\begin{Satz}
\label{aro-theorem}
Let \(A\) be a linear elliptic second-order differential operator defined on a domain \(D\) of \(\R^n\).
Let \(u=(u^1,\ldots,u^n)\) be functions in \(D\) satisfying the inequality
\begin{equation}
\label{aro-voraus}
|Au^j|\leq C\big(\sum_{\alpha,i}\big|\frac{\partial u^i}{\partial x^\alpha}\big|+\sum_i|u^i|\big).
\end{equation}
If \(u=0\) in an open set, then \(u=0\) throughout \(D\).
\end{Satz}

Making use of this result we can prove the following
\begin{Prop}
Let \(\phi\in C^4(M,N)\) be an interpolating sesqui-harmonic map.
If \(\phi\) is harmonic on a connected open set \(W\) of \(M\)
then it is harmonic on the whole connected component of \(M\) which contains \(W\).
\end{Prop}
\begin{proof}
The analytic structure of the interpolating sesqui-harmonic map equation is the
following
\begin{align*}
|\Delta\tau(\phi)|\leq C(|d\phi|^2|\tau(\phi)|+|\tau(\phi)|).
\end{align*}
In order to apply Theorem \ref{aro-theorem} we consider the equation for interpolating sesqui-harmonic maps
in a coordinate chart in the target. The bound on \(|d\phi|^2\) can be obtained by shrinking
the chart if necessary such that \eqref{aro-voraus} holds.
\end{proof}

\subsection{Interpolating sesqui-harmonic maps with vanishing energy-momentum tensor}
In this section we study the qualitative behavior of solutions to \eqref{euler-lagrange}
under the additional assumption that the energy-momentum tensor \eqref{energy-momentum} vanishes
similar to \cite{MR2395125}.
Such an assumption is partially motivated from physics:
In physics one usually also varies the action functional \eqref{energy-functional}
with respect to the metric on the domain and the resulting Euler-Lagrange
equation yields the vanishing of the energy-momentum tensor.

In the following we will often make use of the trace of the energy-momentum tensor \eqref{energy-momentum},
which is given by (where \(n=\dim M\))
\begin{align}
\label{trace-energy-momentum}
\tr T=\delta_1(1-\frac{n}{2})|d\phi|^2+\delta_2\frac{n}{2}|\tau(\phi)|^2
+\delta_2(n-2)\langle d\phi,\nabla\tau(\phi)\rangle.
\end{align}

Note that we do not have to assume that \(\phi\) is a solution of \eqref{euler-lagrange}
in the following.

\begin{Prop}
Let \(\gamma\colon S^1\to N\) be a curve with vanishing energy-momentum tensor.
If \(\delta_1\delta_2>0\) then \(\gamma\) maps to a point.
\end{Prop}
\begin{proof}
Using \eqref{trace-energy-momentum} and integrating over \(S^1\) we find
\begin{align*}
0=\frac{\delta_1}{2}\int_{S^1}|\gamma'|^2ds+\frac{3}{2}\delta_2\int_{S^1}|\tau(\gamma)|^2ds,
\end{align*}
which yields the claim.
\end{proof}

\begin{Prop}
Suppose that \((M,h)\) is a Riemannian surface.
Let \(\phi\colon M\to N\) be a smooth map with 
vanishing energy-momentum tensor. Then \(\phi\) is harmonic.
\end{Prop}
\begin{proof}
Since \(\dim M=2\) we obtain from \eqref{trace-energy-momentum} that \(|\tau(\phi)|^2=0\)
yielding the claim.
\end{proof}

For a higher-dimensional domain we have the following result.

\begin{Prop}
Let \(\phi\colon M\to N\) be a smooth map with 
vanishing energy-momentum tensor. Then the following statements hold:
\begin{enumerate}
 \item If \(\dim M=3\) and \(\delta_1\delta_2<0\) then \(\phi\) is trivial.
 \item If \(\dim M=4\) then \(\phi\) is trivial.
 \item If \(\dim M\geq 5\) and \(\delta_1\delta_2>0\) then \(\phi\) is trivial.
\end{enumerate}

\end{Prop}
\begin{proof}
Integrating \eqref{trace-energy-momentum} we obtain
\begin{align*}
0=\delta_1(1-\frac{n}{2})\int_M|d\phi|^2\dv+\delta_2(2-\frac{n}{2})\int|\tau(\phi)|^2\dv,
\end{align*}
which already yields the result.
\end{proof}

As a next step we rewrite the condition on the vanishing of the energy-momentum tensor.

\begin{Prop}
Let \(\phi\colon M\to N\) be a smooth map and assume that \(n\neq 2\).
Then the vanishing of the energy-momentum tensor is equivalent to
\begin{align}
\label{energy-momentum-vanishing-equivalent}
0=T(X,Y)=&\delta_1\langle d\phi(X),d\phi(Y)\rangle \\
\nonumber&-\delta_2\frac{1}{n-2}|\tau(\phi)|^2h(X,Y) 
-\delta_2\big(\langle d\phi(X),\nabla_Y\tau(\phi)\rangle+\langle d\phi(Y),\nabla_X\tau(\phi)\rangle\big).
\end{align}
\end{Prop}

\begin{proof}
Rewriting the equation for the vanishing of the trace of the energy-momentum tensor \eqref{trace-energy-momentum}
we find
\begin{align*}
\delta_2\langle d\phi,\nabla\tau(\phi)\rangle=\delta_1\frac{(\frac{n}{2}-1)}{(n-2)}|d\phi|^2-\delta_2\frac{n}{2(n-2)}|\tau(\phi)|^2.
\end{align*}
Inserting this into the energy-momentum-tensor \eqref{energy-momentum} yields the claim.
\end{proof}

This allows us to give the following
\begin{Prop}
Let \(\phi\colon M\to N\) be a smooth map with 
vanishing energy-momentum tensor. Suppose that \(\dim M>2\) and \(\operatorname{rank}\phi\leq n-1\).
Then \(\phi\) is harmonic.
\end{Prop}
\begin{proof}
Fix a point \(p\in M\). By assumption \(\operatorname{rank}\phi\leq n-1\) and hence
there exists a vector \(X_p\in\ker d\phi_p\). For \(X=Y=X_p\) we can infer from \eqref{energy-momentum-vanishing-equivalent}
that \(\tau(\phi)=0\) yielding the claim.
\end{proof}

If the domain manifold \(M\) is non-compact we can give the following variant of the 
previous results.

\begin{Prop}
Let \(\phi\colon M\to N\) be a smooth Riemannian immersion with vanishing energy-momentum tensor.
If \(\dim M=2\) then \(\phi\) is harmonic, if \(\dim M=4\) then \(\phi\) is trivial.
\end{Prop}
\begin{proof}
Since \(\phi\) is a Riemannian immersion, we have \(\langle\tau(\phi),d\phi\rangle=0\).
Hence \eqref{trace-energy-momentum} yields
\begin{align*}
0=\delta_1(1-\frac{n}{2})|d\phi|^2+\delta_2(2-\frac{n}{2})|\tau(\phi)|^2,
\end{align*}
which proves the claim.
\end{proof}

\subsection{Conformal construction of interpolating sesqui-harmonic maps}
In \cite{MR1952859} the authors present a powerful construction method 
for biharmonic maps. Instead of trying to directly solve the fourth-order
equation for biharmonic maps they assume the existence of a harmonic map
and then perform a conformal transformation of the metric on the domain
to render this map biharmonic. 
In particular, they call a metric that renders the identity map biharmonic,
a \emph{biharmonic metric}.
In this section we will discuss if the same approach can also be used to
construct interpolating sesqui-harmonic maps.

To this end let \(\phi\colon(M,h)\to (N,g)\) be a smooth map.
If we perform a conformal transformation of the metric on the domain,
that is \(\tilde h=e^{2u}h\) for some smooth function \(u\),
we have the following formula for the transformation of the tension field
\begin{align*}
\tau_{\tilde h}(\phi)=e^{-2u}(\tau_h(\phi)+(n-2)d\phi(\nabla u)),
\end{align*}
where \(\tau_{\tilde h}(\phi)\) denotes the tension field of the map \(\phi\)
with respect to the metric \(\tilde h\). 
In addition, we set \(n:=\dim M\).

Now, suppose that \(\phi\) is a harmonic map with respect to \(h\),
that is \(\tau_h(\phi)=0\), then we obtain the identity
\begin{align*}
\tau_{\tilde h}(\phi)=e^{-2u}(n-2)d\phi(\nabla u).
\end{align*}

This allows us to deduce
\begin{Prop}
Let \(\phi\colon(M,h)\to (N,g)\) be a smooth harmonic map and suppose that \(\dim M\neq 2\).
Let \(\tilde h=e^{2u}h\) be a metric conformal to \(h\).
Then the map \(\phi\colon(M,\tilde h)\to (N,g)\) is interpolating sesqui-harmonic 
if and only if
\begin{align*}
\delta_2\big(\nabla^\ast\nabla d\phi(\nabla u)+(n-6)\nabla_{\nabla u}d\phi(\nabla u)
+2(-\Delta u-(n-4)|du|^2)d\phi(\nabla u) \\
\nonumber+\tr_hR^N(d\phi(\nabla u),d\phi)d\phi\big)
=\delta_1e^{-2u}d\phi(\nabla u).
\end{align*}
\end{Prop}
\begin{proof}
For every \(v\in\Gamma(\phi^\ast TN)\) the following formula holds
\begin{align*}
\nabla^\ast_{\tilde h}\nabla_{\tilde h}v=e^{-2u}(\nabla^\ast_{h}\nabla_{h}v+(n-2)\nabla_{\nabla u}v).
\end{align*}
By a direct calculation we find
\begin{align*}
\nabla_h^\ast\nabla_h(\tau_{\tilde h}(\phi))
=&(n-2)e^{-2u}\big(-2\Delta ud\phi(\nabla u)
+4|du|^2d\phi(\nabla u) \\
&-4\langle\nabla u,\nabla(d\phi(\nabla u))\rangle
+\nabla^\ast\nabla d\phi(\nabla u)\big).
\end{align*}
Together with 
\begin{align*}
\nabla_{\nabla u}\tau_{\tilde h}(\phi)&=(n-2)e^{-2u}(-2|\nabla u|^2d\phi(\nabla u)+\nabla_{\nabla u}d\phi(\nabla u)), \\
R^N(d\phi,\tau_{\tilde h}(\phi))d\phi&=(n-2)e^{-4u}R^N(d\phi,d\phi(\nabla u))d\phi
\end{align*}
this completes the proof.
\end{proof}

In the following we will call a metric that renders the identity map interpolating sesqui-harmonic
an \emph{interpolating sesqui-harmonic metric}.

\begin{Cor}
Let \(\phi\colon(M,h)\to (M,h)\) be the identity map and suppose that \(\dim M\neq 2\).
Let \(\tilde h=e^{2u}h\) be a metric conformal to \(h\).
Then the map \(\phi\colon(M,\tilde h)\to (M,h)\) is interpolating sesqui-harmonic 
if and only if
\begin{align*}
\delta_2(2(-\Delta u-(n-4)|du|^2)\nabla u
+(n-6)\nabla_{\nabla u}\nabla u+\tr_h(\nabla^\ast\nabla)\nabla u
+\operatorname{Ric}^M(\nabla u))
=\delta_1e^{-2u}\nabla u.
\end{align*}
\end{Cor}
We now rewrite this as an equation for \(\nabla u\).

\begin{Prop}
Let \(\phi\colon(M,h)\to (M,h)\) be the identity map and suppose that \(\dim M\neq 2\).
Let \(\tilde h=e^{2u}h\) be a metric conformal to \(h\) and set \(\nabla u=\beta\).
Then the map \(\phi\colon(M,\tilde h)\to (M,h)\) is interpolating sesqui-harmonic 
if and only if
\begin{align}
\label{equation-metric-beta}
-\Delta\beta=2(d^\ast\beta-(n-4)|\beta|^2)\beta+(n-6)\frac{1}{2}d|\beta|^2
+2\operatorname{Ric}^M(\beta^\sharp)^\flat-\frac{\delta_1}{\delta_2}e^{-2u}\beta,
\end{align}
where \(-\Delta=dd^\ast+d^\ast d\) is the Laplacian acting on one-forms.
\end{Prop}
\begin{proof}
The Laplacian acting on one-forms satisfies the following Weitzenböck identity
\begin{align*}
\Delta\beta(X)=\tr_h(\nabla^2)\beta(X)-\beta\operatorname{Ric}(X),
\end{align*}
where \(X\) is a vector field.
In addition, note that (see \cite[Proof of Proposition 2.2]{MR1952859} for more details)
\begin{align*}
\nabla_{\nabla u}\nabla u=\frac{1}{2}d|\beta|^2,
\end{align*}
which completes the proof.
\end{proof}

\begin{Prop}
Let \((M,h)\) be a compact manifold of strictly negative Ricci curvature with \(\dim M>2\) 
and assume that \(\delta_1\delta_2>0\). Then there does not exist an interpolating sesqui-harmonic metric 
that is conformally related to \(h\) except a constant multiple of \(h\).
\end{Prop}
\begin{proof}
Note that our sign convention for the Laplacian is different from the one used in \cite{MR1952859}.
We define the one-form \(\theta:=e^{-2u}\beta\). By a direct calculation we then find using \eqref{equation-metric-beta}
\begin{align*}
\Delta\theta=-\frac{1}{2}(n-2)e^{2u}d|\theta|^2-2\operatorname{Ric}^M(\theta^\sharp)^\flat+\frac{\delta_1}{\delta_2}e^{-2u}\theta.
\end{align*}
In addition, we have
\begin{align*}
\Delta\frac{1}{2}|\theta|^2&=\langle\Delta\theta,\theta\rangle+|\nabla\theta|^2+\operatorname{Ric}(\theta^\sharp,\theta^\sharp) \\
&=-\frac{1}{2}(n-2)e^{2u}\langle d|\theta|^2,\theta\rangle+|\nabla\theta|^2-\operatorname{Ric}(\theta^\sharp,\theta^\sharp)+\frac{\delta_1}{\delta_2}e^{-2u}|\theta|^2.
\end{align*}
The claim then follows by the maximum principle.
\end{proof}

\begin{Bem}
In contrast to the case of biharmonic maps \eqref{equation-metric-beta} contains also
a term involving \(u\) on the right hand side. This reflects the fact that both
harmonic and biharmonic maps on its own have a nice behavior under conformal deformations
of the domain metric, whereas interpolating sesqui-harmonic maps do not.
This prevents us from making a connection between interpolating sesqui-harmonic metrics
and isoparametric functions as was done in \cite{MR1952859}
for biharmonic maps.
\end{Bem}

\subsection{A Liouville-type theorem for interpolating sesqui-harmonic maps between complete manifolds}
In this section we will prove a Liouville-type theorem
for solutions of \eqref{euler-lagrange} between complete Riemannian manifolds
generalizing a similar result for biharmonic maps from \cite{MR3205803}.
For more Liouville-type theorems for biharmonic maps see \cite{MR3834926,branding2018nonexistence} and references therein.

To this end we will make use of the following result due to Gaffney \cite{MR0062490}:
\begin{Satz}
\label{gaffney}
Let \((M,h)\) be a complete Riemannian manifold. If a \(C^1\) one-form \(\omega\)
satisfies 
\begin{align*}
\int_M|\omega|\dv<\infty \qquad\text{ and }\qquad \int_M|\delta\omega|\dv<\infty
\end{align*}
or, equivalently, a \(C^1\) vector field \(X\) defined by \(\omega(Y)=h(X,Y)\),
satisfies
\begin{align*}
\int_M|X|\dv<\infty \qquad\text{ and }\qquad \int_M\operatorname{div}(X)\dv<\infty,
\end{align*}
then
\begin{align*}
\int_M(\delta\omega)\dv=\int_M\operatorname{div}(X)\dv=0.
\end{align*}
\end{Satz}

\begin{Satz}
Let \((M,h)\) be a complete non-compact Riemannian manifold and \((N,g)\) a manifold with non-positive sectional curvature.
Let \(\phi\colon M\to N\) be a smooth solution of \eqref{euler-lagrange} and \(p\) be a real constant satisfying
\(2\leq p<\infty\).
\begin{enumerate}
 \item If \(\delta_1\delta_2>0\) and \begin{align*}
        \int_M|\tau(\phi)|^p\dv<\infty,\qquad \int_M|d\phi|^2\dv<\infty
       \end{align*}
then \(\phi\) must be harmonic.
\item If \(\delta_1\delta_2>0\), \(\vol(M,h)=\infty\) and
\begin{align*}
        \int_M|\tau(\phi)|^p\dv<\infty
\end{align*}
then \(\phi\) must be harmonic.
\end{enumerate}
\end{Satz}

\begin{proof}
We choose a cutoff function  \(0\leq\eta\leq 1\) on \(M\) that satisfies
\begin{align*}
\eta(x)=1\textrm{ for } x\in B_R(x_0),\qquad \eta(x)=0\textrm{ for } x\in B_{2R}(x_0),\qquad |\nabla\eta|\leq\frac{C}{R}\textrm{ for } x\in M,
\end{align*}
where \(B_R(x_0)\) denotes the geodesic ball around the point \(x_0\) with radius \(R\).

We test the interpolating sesqui-harmonic map equation \eqref{euler-lagrange} with \(\eta^2\tau(\phi)|\tau(\phi)|^{p-2}\) and find
\begin{align*}
\eta^2|\tau(\phi)|^{p-2}\langle\Delta\tau(\phi),\tau(\phi)\rangle=
\eta^2|\tau(\phi)|^{p-2}\langle R^N(d\phi(e_\alpha),\tau(\phi))d\phi(e_\alpha),\tau(\phi)\rangle
+\frac{\delta_1}{\delta_2}\eta^2|\tau(\phi)|^{p} 
\geq 0,
\end{align*}
where we made use of the assumptions on the curvature of the target and the signs of \(\delta_1\) and \(\delta_2\).
Integrating over \(M\) and using integration by parts we obtain
\begin{align*}
\int_M\eta^2|\tau(\phi)|^{p-2}\langle\Delta\tau(\phi),\tau(\phi)\rangle\dv=&-2\int_M\langle\nabla\tau(\phi),\tau(\phi)\rangle|\tau(\phi)|^{p-2}\eta\nabla\eta\dv \\
&-(p-2)\int_M\eta^2|\langle\nabla\tau(\phi),\tau(\phi)\rangle|^2|\tau(\phi)|^{p-4}\dv \\
&-\int_M\eta^2|\nabla\tau(\phi)|^2|\tau(\phi)|^{p-2}\dv
\\
\leq&\frac{C}{R^2}\int_M|\tau(\phi)|^{p}\dv
-\frac{1}{2}\int_M\eta^2|\nabla\tau(\phi)|^2|\tau(\phi)|^{p-2}\dv\\
&-(p-2)\int_M\eta^2|\langle\nabla\tau(\phi),\tau(\phi)\rangle|^2|\tau(\phi)|^{p-4}\dv,
\end{align*}
where we used Young's inequality and the properties of the cutoff function \(\eta\).
Combining both equations we find
\begin{align*}
\frac{1}{2}\int_M\eta^2|\nabla\tau(\phi)|^2|\tau(\phi)|^{p-2}\dv\leq &\frac{C}{R^2}\int_M|\tau(\phi)|^{p}\dv
-(p-2)\int_M\eta^2|\langle\nabla\tau(\phi),\tau(\phi)\rangle|^2|\tau(\phi)|^{p}\dv.
\end{align*}
Letting \(R\to\infty\) and using the finiteness assumption of the \(L^p\) norm of the tension field we 
may deduce that \(\tau(\phi)\) is parallel and thus has constant norm.

To establish the first claim of the theorem we make use of Theorem \ref{gaffney}.
We define a one-form \(\omega\) by
\begin{align*}
\omega(X):=|\tau(\phi)|^{\frac{p}{2}-1}\langle d\phi(X),\tau(\phi)\rangle,
\end{align*}
where \(X\) is a vector field on \(M\). 
Note that
\begin{align*}
\int_M|\omega|\dv\leq\int_M|d\phi||\tau(\phi)|^{\frac{p}{2}}\dv\leq\big(\int_M|d\phi|^2\dv\big)^\frac{1}{2}\big(\int_M|\tau(\phi)|^p\dv\big)^\frac{1}{2}<\infty.
\end{align*}
Using that \(|\tau(\phi)|\) has constant norm we find by a direct calculation that \(\delta\omega=|\tau(\phi)|^{\frac{p}{2}+1}\).
Again, since \(|\tau(\phi)|\) has constant norm and the \(L^p\)-norm of \(\tau(\phi)\) is bounded, we find that \(|\delta\omega|\)
is integrable over \(M\). By application of Theorem \ref{gaffney} we can then deduce that \(\tau(\phi)=0\).

To prove the second claim, we note that \(\vol(M,h)=\infty\) and \(|\tau(\phi)|\neq 0\) give
\begin{align*}
\int_M|\tau(\phi)|^p\dv=|\tau(\phi)|^p\vol(M,h)=\infty,
\end{align*}
yielding a contradiction.
\end{proof}

\par\medskip
\textbf{Acknowledgements:}
The author gratefully acknowledges the support of the Austrian Science Fund (FWF) 
through the project P30749-N35 ``Geometric variational problems from string theory''.
\bibliographystyle{plain}
\bibliography{mybib}
\end{document}